\documentclass[12pt]{article}

\usepackage{latexsym,amssymb}

\def\tr{^{\rm T}}
\def\Real{\mathbb R}

\def\dst{\displaystyle}

\def\qmx#1{\left(\matrix{#1}\right)}

\def\beeq#1{\begin{equation}{#1}\end{equation}}

\def\ba{\begin{array}}
\def\ea{\end{array}}
\def\eqa{\begin{eqnarray}}
\def\eqe{\end{eqnarray}}

\newtheorem{proposition}{Proposition}

\newtheorem{lemma}{Lemma}

\newenvironment{proof}{\medskip\noindent{\it Proof. }}{ \medskip}
\newenvironment{remark}{\medskip\noindent{\it Remark. }}{
\medskip}
\newenvironment{example}{\medskip\noindent{\it Example. }}{
\medskip}

\begin{document}

\title{Adaptive Observers as Nonlinear Internal Models\thanks{This work
was supported in part by ONR under grant
 N00014-03-1-0314, and by MIUR. Corresponding
 author: Dr. Lorenzo Marconi, Tel: 0039.051.2093788, Fax:
 0039.051.2093073,
 email: lmarconi@deis.unibo.it.
 }}

\author{F. Delli Priscoli $^{\dag}$, L. Marconi $^{\circ}$,
A.Isidori $^{\dag \ddag \circ}$}

\date{}

\maketitle

\small {

\begin{center}
$^{\dag}$Dipartimento di Informatica e Sistemistica,
 Universit\`{a} di Roma ``La Sapienza'', \\00184 Rome, ITALY.

$^{\circ}$ C.A.SY. -- Dipartimento di Elettronica, Informatica e
Sistemistica, University of Bologna, \\40136 Bologna, ITALY.

$^{\ddag}$Department of Electrical and Systems Engineering,
Washington University, \\St. Louis, MO 63130.
\end{center}
}
\smallskip
\normalsize

\begin{abstract}
This paper shows how the theory of nonlinear adaptive observers
can be effectively used in the design  of internal models for
nonlinear output regulation. The theory substantially enhances the
existing results in the context of {\em adaptive} output
regulation, by allowing for not necessarily stable zero dynamics
of the controlled plant and by weakening  the standard assumption
of having the steady state control input generated by a linear
system.
\end{abstract}

{\em keywords}: Adaptive Observers, Internal Model, Regulation,
Tracking, Nonlinear Control.

\section{Introduction}\label{sec1}

The problem of controlling the output of a system so as to achieve
asymptotic tracking of prescribed trajectories and/or asymptotic
rejection of disturbances, sometimes known as the {\em
servomechanism problem}, is continuing to attract a good deal of
interest in control theory. In the last decade or so, most of the
efforts have been addressed toward the design of controllers which
solve this problem in the case of plants modelled by nonlinear
differential equations. Viewed as a nonlinear design problem, some
of the original features (such as, for instance, the conservation
of the desired steady-state features in spite of plant parameter
variations, otherwise known as ``robustness" property, and the
necessity of an ``internal model" in any robust regulator) tend to
lose their specific connotation. Rather, they merge with other
relevant issues in feedback design for nonlinear systems, notably
the guarantee of convergence (for certain state variables) or
boundedness (for other state variables) once a fixed set of
initial data is given. After all, a problem of steering certain
variables to a desired target value in the presence of exogenous
stimuli generated by an autonomous ``exosystem", in a nonlinear
context, can be viewed as a problem of adaptive control.
Classically, adaptation is sought with respect to uncertain but
constant parameters (which, of course, can be seen as exogenous
inputs generated by a ``trivial" autonomous exosystem) but if the
parameters in question are obeying some fixed differential
equation, it would not be inappropriate to continue to call a
problem of this kind a (generalized) problem of adaptive control.
It is because of this observation that an increasing use of
methods and techniques from (nonlinear) adaptive control should be
expected, so long as newer results in this research area will be
generated.

A specific feature of the problem in question is that, to achieve
the desired steady-state features (perfect tracking), the
controller should be able to generate a family of special inputs
(those which, in fact, secure perfect tracking). If the external
stimuli are constant (as in the case of uncertain constant plant
parameters), this generator is simply provided by a bank of
integrators (in the case of adaptive control, the state of each of
such integrator is a {\em parameter estimate}). These integrators
are to be ``controlled" by appropriate feedback laws, so as to
achieve the desired convergence and/or boundedness properties (in
the case of adaptive control, this is the design of the
``adaptation laws"). In a general servomechanism problem, the
setting is very much the same: a model which generates all inputs
needed to obtain perfect tracking is first found (the ``internal
model") and then some (generalized) stabilization law is
superimposed to complete the design. In this respect, the design
of a controller that solves the servomechanism problem is split in
two parts: the design of an internal model and the design of a
stabilizer. It must be observed, though, that the design of the
former must be done with an eye to our ability to find the latter.
This is a fact that was well understood in adaptive control. The
addition of appropriate dynamics (notably the so-called ``filtered
transformations") are not strictly speaking needed to provide the
dynamics of the unknown parameters (a trivial dynamics in that
case), but rather are introduced to address the issue of
stability.

This being said, it is natural to expect an increasing interaction
between the work on nonlinear adaptive control, nonlinear
stabilization and nonlinear servomechanism theory. This
interaction has already manifested itself in a number of recent
contributions (such as \cite{SIM}, \cite{BI03}, \cite{JMHH},
\cite{CH2}, \cite{BIP}, \cite{PWN},  whose results cannot be
reviewed here for obvious reasons) and will manifest more in the
contributions to come. In this paper, we wish to propose a method
for design of an internal model which is based on some classical
results in adaptive control: the design of adaptive observers for
nonlinear systems which are linearizable by output injection. This
method enables us to solve a servomechanism problem when the
inputs needed to achieve perfect tracking can be seen as generated
by a nonlinear system, linearizable by output injection, with
possibly unknown coefficients. Allowing for unknown coefficients
in this model automatically settles the issue of uncertain plant
parameters (the classical ``robustness" issue) as well the issue
of parameter uncertainties in the exosystem (an outstanding design
problem first addressed and solved, for a special class of
systems, in \cite{SIM}). Although the inspiration for the design
is taken from an existing result in adaptive control, the
application to the specific context of servomechanism theory looks
pretty new and worth being pursued.

\section{Class of systems and main assumptions}\label{SecPA}

\subsection{Preliminaries}

The purpose of this paper is to show how the theory of nonlinear
adaptive observed can be effectively used in the design of
adaptive output regulators for nonlinear systems. To motivate why
and how adaptive observers play an important role in this design
problem, it suffices to address the simplified case in which the
controlled plant has relative degree 1 between the {\em control
input} and the {\em regulated output}. This is what is done here,
for reason of space. The extension of the design methodology to
system having higher relative degree can be found in the more
extended paper \cite{DMISIAM}, along with some additional
technical details.

 Consider a system  modelled by equations of the form
 \beeq{\label{pla1}\ba{rcl} \dot z &=& f_0(\varrho,w,z) +
 f_1(\varrho,w, z,e)e\\
 \dot e &=& q(\varrho,w,z,e) + u\,,
 \ea}
 with state $(z,e)\in \Real^n \times \Real$,
control input $u\in \Real$, regulated output $e\in \Real$,
 in which the exogenous
inputs $\varrho\in \Real^p$ and $w\in \Real^s$ are generated by an
exosystem modelled by equations of the form \beeq{\label{exo1}
 \ba{rcl}
 \dot \varrho&=&0\\
 \dot w &=& s(\varrho,w)\,.
 \ea
 }
In this model, $\varrho \in \Real^p$ is a vector of constant
uncertain parameters, the aggregate of a finite set of uncertain
parameters affecting the controlled plant and another, possibly
different, set of uncertain parameters affecting the generator of
the exogenous input $w$. Note that system (\ref{pla1}) has
relative degree 1 between control input $u$ and regulated output
$e$.

The functions $f_0(\cdot), f_1(\cdot), q(\cdot), s(\cdot)$ in
(\ref{pla1}) and (\ref{exo1}) are assumed to be at least
continuously differentiable. The initial conditions of
(\ref{pla1}) range on a set $Z\times E$, in which $Z$ is a fixed
{\em compact} subset of $\Real^n$ and $E=\{e\in \Real:|e|\le c\}$,
with $c$ a fixed number. The initial conditions of the exosystem
(\ref{exo1}) range on set $P\times W$ in which $P$ and $W$ are
{\em compact} subsets of $\Real^p$ and, respectively, $\Real^s $.
In this framework the problem of output regulation is to design an
output feedback regulator of the form
 \[
 \ba{rcl}
 \dot \zeta &=& \varphi(\zeta, e)\\
 u &=& \gamma(\zeta,e)
 \ea
 \]
 such that {\em for all initial conditions $(\varrho(0),w(0)) \in
 P\times W$ and $(z(0),e(0)) \in Z \times E$ the
 trajectories of the closed-loop system are bounded and $\lim_{t \rightarrow \infty}
 e(t)=0$.}

We retain in this paper some ideas introduced in \cite{BI03}, to
which -- to avoid duplications -- the reader is referred. Among
the concepts introduced and/or summarized in that paper, the
notion of {\em omega limit set $\omega({\bf S})$ of a set ${\bf
S}$} plays a major role. This concept is a deep generalization of
the classical concept, due to Birkhoff, of omega limit set of a
point and provides a rigorous definition of steady-state response
in a nonlinear system (see \cite{BI03} for details).

\begin{remark} The regulated variable $e$ of (\ref{pla1}) may coincide with
 the physical ``controlled" output of a given
plant, or may as well represents a {\em tracking error}, namely
the difference between a physical ``controlled" output and its
``reference" behavior. Thus, the problem under consideration
includes problems of tracking as well as problems of disturbance
attenuation. $\triangleleft$
\end{remark}

 \subsection{Basic hypotheses}

Augmenting (\ref{pla1}) with (\ref{exo1}) yields a system which,
viewing $u$ as input and $e$ as output, has relative degree $1$.
The associated ``augmented" zero dynamics, which is forced by the
control
 \beeq{ \label{friend1}  c(\varrho,
 w, z) = -q(\varrho,w,z,0)\,,
 }
is given by \beeq{\label{zerodynpla}\ba{rcl}
 \dot \varrho &=&0\\[2mm]
 \dot w &=& s(\varrho,w)\\[2mm]
 \dot z &=& f_0(\varrho,w,z) \,.\ea
 }

Occasionally, throughout the paper, we will rewrite the latter as
 \beeq{\label{zerobold}
 \dot {\bf z} = {\bf f}_0({\bf z})\,,}
 where
 ${\bf z}={\rm col}(\varrho,w,z)$.
 Accordingly, we set ${\bf Z}= P\times W\times Z$ and
 replace $c(\varrho,
 w, z)$ by $c({\bf z})$ in (\ref{friend1}).

In what follows, we retain three of the basic assumptions that
were introduced in \cite{BI03} and express certain properties of
the augmented zero dynamics (\ref{zerodynpla}). The assumptions in
question are the following ones:

\medskip\noindent {\em Assumption (i)}\,: the set $P\times W$ is a differential submanifold (with
boundary) of $\Real^{p} \times \Real^s$, invariant for
(\ref{exo1}). $\triangleleft$

\medskip\noindent {\em Assumption (ii)}\,:
there exists a compact subset $\mathcal Z$ of $P\times W \times
\Real^n$ which contains the positive orbit of the set ${\bf Z}$
under the flow of (\ref{zerobold}), and $\omega({\bf Z})$ is a
differential submanifold (with boundary) of $P\times W \times
\Real^n$. Moreover there exists a number $d_1>0$ such that
 \[
 {\bf z} \in
P\times W \times \Real^n\,, \quad \mbox{dist}({\bf z}, \omega({\bf
Z})) \leq d_1 \qquad \Rightarrow \qquad {\bf z} \in {\bf Z}\,.
\quad \triangleleft
\]

\begin{remark} Since the
positive orbit of the set $\bf Z$ under the flow of
(\ref{zerobold}) is bounded by hypothesis, the set $\omega({\bf
Z})$ is a nonempty, compact and invariant subset of $P\times W
\times \Real^n$ which uniformly attracts all trajectories of
(\ref{zerobold}) with initial conditions in ${\bf Z}$. It can also
be shown (as in \cite{BI03}) that  for every $(\varrho,w) \in
P\times W$ there is $z \in \Real^n$ such that $(\varrho,w, z) \in
\omega({\bf Z})$. $\triangleleft$
\end{remark}

In what follows, for convenience, the set $\omega({\bf Z})$ will
be simply denoted as ${\mathcal A}_0$. The last condition in
assumption (ii) implies that ${\mathcal A}_0$ is stable in the
sense of Lyapunov. The next hypothesis is that the set ${\mathcal
A}_0$ is locally exponentially attractive.

\medskip\noindent {\em Assumption (iii)}\,: There exist
$M \geq 1$, $a>0$ and $d_2 \leq d_1$ such that
\[
{\bf z}_0 \in P\times W \times \Real^n\,, \quad \mbox{dist}({\bf
z}_0, {\mathcal A}_0) \leq d_2 \qquad \Rightarrow \qquad
 \mbox{dist}({\bf z}(t,{\bf z}_0), \, {\mathcal A}_0) \leq M e^{- a t}
 \mbox{dist}({\bf z}_0 , \, {\mathcal A}_0)
\]
in which ${\bf z}(t,{\bf z}_0)$ denotes the solution of
(\ref{zerobold}) passing through ${\bf z}_0$ at time $t=0$.
$\triangleleft$

\medskip
As is well known, a fundamental step in the construction of a
controller which solves the problem of output regulation is the
design on an ``internal model", i.e. an autonomous dynamical
system which generates all ``feed-forward inputs capable to secure
perfect tracking" (in the present case, the set of all inputs of
the form $u(t) = c({\bf z}(t))$, with ${\bf z}(t)$ a trajectory of
the restriction of (\ref{zerobold}) to ${\mathcal A}_0$). In this
respect, a relevant contribution to the design of internal models
has been the methodology originally developed by Huang and
co-authors (see e.g. \cite{HL,JH1,JH2}). As a matter of fact, in
these works  a precise characterization is provided of when the
function $u(t)$ in question satisfies a linear differential
equation
\[
u^{(d)}(t)+a_{d-1}u^{(d-1)}(t)+ \cdots + a_{1}u^{(1)}(t)+a_0u(t)
=0
\]
which can then be taken as internal model and provides the
fundamental core for the design of regulators. It was only
recently, though, that attempts to weaken this condition have been
taken. In \cite{SIM}, for instance, the case in which the
coefficients of the above equation may depend on the uncertain
(constant) parameter $\varrho$ was considered. Alternative and
relevant extensions have also been proposed in the recent works
\cite{CH1,CH2}. The extension to the case in which the equation
above is nonlinear (but $\varrho$-independent) was treated in
\cite{BI03bis}.

The present paper proposes an extension, in the design of internal
models, to the case in which the autonomous dynamical system that
generates $u(t)$ is nonlinear but ``linearizable by output
injection", with coefficients which possibly depend on the
uncertain (constant) parameter $\varrho$. This case, to the best
of our knowledge, has not been treated before by other authors.
The assumption which characterizes when the extension in question
is possible is the following one.

\medskip\noindent {\em Assumption (iv)}\,:
 there exist a positive integer $d$, a $C^1$ map
 \[
 \ba{lrcl}
 \tau \quad : \quad &  {\mathcal Z} &\rightarrow& \Real^d\\
              & {\bf z} &\mapsto& \tau({\bf z})\,,
 \ea
 \]
 a $C^0$ map
 \[
 \ba{lrcl}
 \theta \quad : \quad &  {P} &\rightarrow& \Real^q\\
              & \varrho &\mapsto& \theta(\varrho)\,,
 \ea
 \]
 an observable pair $(A,C) \in \Real^{d \times d} \times
 \Real^{1 \times d}$, and two $C^1$ maps $\phi: \Real \to \Real^d$ and
 $\Omega: \Real \to \Real^{d \times q}$ such that the following
 identities (which we call
 {\em immersion property})
 \beeq{ \label{adaptimmers1}
 {\partial \tau \over \partial {\bf z}} \, {\bf f}_0({\bf z})
   =A\, \tau({\bf z}) + \phi(C\tau({\bf z}))
  + \Omega(C\tau({\bf z}))\, \theta(\varrho)
 }
 \beeq{\label{adaptimmers0}
 c({\bf z}) = C\, \tau({\bf z})\,
 }
 hold for all ${\bf z}\in {\mathcal A}_0$, $\varrho \in P$. $\triangleleft$

 \begin{remark}
 Without loss of generality (see \cite[page 208]{MaToBook}),
 we can assume throughout that the matrices $A$ and $C$ in
 (\ref{adaptim}) have the form
 \[
 A =\qmx{0 & 1 & 0 &\cdots &0\cr 0 & 0 & 1 &\cdots &0\cr \cdot&
 \cdot & \cdot & \cdots &\cdot \cr 0 & 0 & 0 &\cdots &1\cr 0 & 0 &
 0& \cdots& 0\cr},\qquad C = \qmx{1 & 0 & 0 &\cdots &0\cr}.
 \]
 Furthermore, note that since the maps
 $\Omega(\cdot)$ and $\phi(\cdot)$ are continuously differentiable and the relations
 (\ref{adaptimmers1}) -- (\ref{adaptimmers0})
 are supposed to hold over the compact set ${\mathcal A}_0$, it can be assumed without loss of
 generality that functions $\phi(\cdot)$ and $\Omega(\cdot)$ have compact support.
 This being the case, the functions in question can be assumed
 {\em globally Lipschitz}, i.e. there exist $L_\phi$ and $L_\Omega$ such that
 \[
 |\phi(s_1)-\phi(s_2)| \le L_\phi|s_1-s_2|, \qquad
 |\Omega(s_1)-\Omega(s_2)| \le L_\Omega|s_1-s_2|,
 \]
 for all $s_1, s_2$. $\triangleleft$
 \end{remark}

Assumption (iv) can be rephrased by saying
 that for each initial condition ${\bf z}(0) \in {\mathcal A}_0$ of (\ref{zerobold}),
 there is a pair
 $\xi(0),\theta$ such that the control input
 $u(t)=c({\bf z(t)})$  can be seen as output of a system
 of the form
 \beeq{ \label{adaptim}
 \ba{rcl}\dot \xi &=& A\xi + \phi(y) +
 \Omega(y)\theta\\[2mm]\dot \theta &=& 0\\[2mm]
 y&=&C\xi\,. \ea
 }

\subsection{Comparison with earlier work on adaptive regulation}

It may be useful to compare the assumption above with some of the
standing assumptions in the earlier work \cite{SIM} on adaptive
regulation. In framework of \cite{SIM}, the lower subsystem of
(\ref{exo1}) was assumed to be a neutrally stable linear system
\[
 \dot w = S(\varrho)w\,.
\]

Assumption III.1 of \cite{SIM} was the existence of a map $\zeta:
P\times W \to \Real^n$ satisfying
\[
{\partial \zeta \over \partial w}S(\varrho)w =
f_0(\varrho,w,\zeta(\varrho,w))\,,\] or, what is the same, of a
map $\zeta: P\times W \to \Real^n$ whose graph is a compact
invariant set for (\ref{zerodynpla}). Assumption V.1 of \cite{SIM}
was that the graph of $\zeta$ is globally asymptotically and
locally exponentially stable. As a consequence, it is readily seen
that Assumptions III.1 and V.1 of \cite{SIM} imply the fulfillment
of assumptions (ii) and (iii) above and, in particular, the fact
that the set ${\mathcal A}_0$ is the graph of the map $\zeta$.

Assumption III.2 of \cite{SIM} was the existence of an integer
$q$, of a $q\times q$ matrix $\Phi(\varrho)$, continuously
depending on $\varrho$, of a $1 \times q$ matrix $\Gamma$ such
that the pair $(\Phi(\varrho), \Gamma)$ is observable for all
$\varrho$, and of a $C^1$ map $\tau: P\times W \to \Real$
satisfying the pair of conditions
\[\ba{rcl}
\dst {\partial \tau \over \partial w}S(\varrho)w &=&
\Phi(\varrho)\tau(\rho,w)\\[2mm]
c(\varrho,w,\zeta(\varrho,w)) &=& \Gamma\tau(\rho,w)\,,\ea\] which
were referred to as conditions of ``immersion" into a {\em linear}
(observable) system. It is shown now that this assumption
automatically yields assumption (iv) above.

Since the pair $(\Phi(\varrho), \Gamma)$ is observable for all
$\varrho$, there exist a matrix $M(\varrho)$, continuously
depending on $\varrho$, and continuous functions
$\theta_1(\varrho), \ldots, \theta_q(\varrho)$ such that
\[
 A(\varrho) := M(\varrho) \Phi(\varrho)M(\varrho)^{-1} = \left ( \ba{ccccc}
 \theta_1(\varrho) & 1 & 0 &\ldots & 0\\
 \theta_2(\varrho) & 0 & 1 &\ldots & 0\\
 \vdots & \vdots & \vdots & \ddots &  \vdots\\
 \theta_q(\varrho) & 0 & 0 & \ldots & 1
 \ea \right )\]
 \[ C:=\Gamma M(\varrho)^{-1}=\left ( \ba{ccccc} 1 & 0 & 0 & \ldots & 0
 \ea \right )\,.
 \]
 The matrix $A(\varrho)$ can be written in the form $A(\varrho) = A + \theta(\varrho) C$
 in which $\theta(\varrho)$ is the column vector defined by the first column of $A(\varrho)$
 and the pair $(A,C)$ is observable. Setting $\bar \tau(\varrho,w)
 = M(\varrho)\tau(\varrho,w)$, and $\Omega(x) = x I_q$, it is
 readily seen that
 \[\ba{rcl}
\dst {\partial \bar \tau \over \partial w}S(\varrho)w &=&
A\bar \tau(\rho,w)+\Omega(C\bar \tau(\rho,w))\theta(\varrho)\\[2mm]
c(\varrho,w,\zeta(\varrho,w)) &=& C\bar \tau(\rho,w)\,.\ea\] This,
bearing in mind the fact that ${\mathcal A}_0$ is the graph of
$\zeta$, shows that assumption III.2 of \cite{SIM} implies the
fulfillment of assumption (iv) above.

The main limitation of the approach of \cite{SIM} was the
assumption of ``immersion" into a {\em linear} system, namely the
assumption that the set of  inputs capable to secure perfect
tracking is a subset of the set of solutions of a suitable {\em
linear} (though $\varrho$-dependent) differential equation. To
show that this limitation is overcome by the approach presented in
this paper, we give here a simple but significant example.

\begin{example}
Consider the regulation problem
\[\ba{rcl}
\dot w_1 &=& w_2\\
\dot w_2 &=& - \omega^2 w_1\\
\dot z_1 &=& z_2\\
\dot z_2 &=& -\sigma z_1-(z_1^2-1)z_2 -w_1+e\\
\dot e &=& -\mu z_1 + u \ea
\]
in which $\omega$, $\sigma$ and $\mu$ are uncertain parameters
ranging over compact sets. Note that $\omega$ is the uncertain
frequency of a linear exosystem, while $\sigma$ characterizes the
frequency of a stable limit cycle embedded in the zero dynamics of
the controlled plant. Note also that the equilibrium $z=0$ of the
zero dynamics is unstable.

Set $\varrho = {\rm col}(\omega, \sigma, \mu)$ and ${\bf z}={\rm
col}(\varrho,w,z)$. The function $c({\bf z})$ in this case is the
function $c({\bf z})=\mu z_1$. Set now
\[
d(x,\mu) = {1 \over \mu^2}x^2 - 1\, \] and \[ \ba{rcl}\tau_1({\bf
z})&=& \mu z_1\\ \tau_2({\bf z})&=& \mu z_2 +  \dst \int_0^{\mu
z_1}d(x,\mu)dx\\ \tau_3({\bf z})&=& \mu w_1 - \omega^2\mu z_1\\
\tau_4({\bf z})&=& \mu w_2 - \omega^2\Bigl[\mu z_2 + \dst
\int_0^{\mu z_1}d(x,\mu)dx\Bigr]\,.\ea\]

A simple calculation shows that the map  $\tau({\bf z}) $ thus
defined satisfies conditions (\ref{adaptimmers1}) and
(\ref{adaptimmers0}), with $(A,C)$ an observable pair and
\[\phi_0(y)= \qmx{y \cr 0 \cr
0 \cr 0 \cr }, \qquad \Omega(y)= \qmx{-y^3 & 0 & 0 & 0\cr 0 & -y &
0 & 0 \cr 0& 0& -y &0\cr 0&0& 0& -y}, \qquad
\theta(\varrho)=\qmx{1/\mu \cr \sigma \cr -\omega^2+\sigma \cr
-\omega^4 - \omega^2\sigma\cr}.\]
\end{example}

\section{ The adaptive internal model}

\subsection{The structure of the regulator}

As pointed out in \cite{BI03bis}, the design of internal models
can be reduced to the design of observers ({\em nonlinear}
observers, in that case). In the presence of uncertain parameters
in the exosystem, an {\em adaptive} observer should be used.
Bearing in mind the construction of adaptive observers for systems
which are linearizable by output injection pioneered by Bastin and
Gevers \cite{BaGe}, we consider in what follows a candidate
controller of the form
 \beeq{\label{control}\ba{rcl} u &=& \xi_1 + v\\[2mm]
 \dot \xi &=& A\xi + \phi(\xi_1) + \Omega(\xi_1)\hat \theta
 +H(X,\xi_1)v - M(X) \mbox{dzv}_\ell(\hat \theta)\\[2mm]
 \dot {\hat \theta} &=& \beta(X,\xi_1) v -\mbox{dzv}_\ell(\hat \theta)\\[2mm]
 \dot X &=& FX + G\Omega(\xi_1)\ea
 }
 in which $\xi_1$ denotes the first component of $\xi$, the matrix
$X$ is a $(d-1)\times q$ matrix, $M(X)$ is a $d\times q$ matrix
defined as
 \[
 M(X) = \qmx{0 \cr X\cr}\,,
 \]
 while the vectors $H(X,\xi_1),
\beta(X,\xi_1)$ and the matrices $F, G$ have the form described
below.  The vector-valued {\em dead-zone}  function
dzv$_\ell(\cdot)$ is defined as
\[
    \mbox{dzv}_\ell(\mbox{col}(s_1,\ldots,s_q)) = \mbox{col}(\mbox{dz}_\ell(s_1), \ldots,
    \mbox{dz}_\ell(s_q))
\]
in which dz$_\ell(\cdot)$ is any continuously differentiable
function satisfying
 \beeq{\label{deftas}
\mbox{dz}_\ell(x) =
 \left \{ \ba{ll} 0 & \quad \mbox{if } |x|
 \leq \ell\\
 x & \quad \mbox{if } |x| \geq \ell +1
  \ea \right .
 }
 and the amplitude $\ell$ of the dead-zone is chosen so that
\[
\ell > \max_{\varrho \in P} |\theta(\varrho)|\,.
\]
 This controller can be viewed as a ``copy" of
(\ref{adaptim}), corrected by an ``innovation term", augmented
with an ``adaptation law" for ${\hat \theta}$ and with a ``filter"
which generates the ``auxiliary state" $X$. The additional input
$v$, which is a ``stabilizing control", will eventually be taken
as $v = -k e$.

The choice of the functions $H(X,\xi_1), \beta(X,\xi_1)$ and the
matrices $F, G$ of (\ref{control}) is inspired by certain
calculations used in \cite{BaGe} and \cite{MaTo92} for the design
of adaptive nonlinear observers. Define new variables
\beeq{\label{newvar}\ba{rcl}
\tilde \theta &=& \hat \theta - \theta(\varrho)\\[2mm] \eta &=& \xi -M(X)\tilde
\theta\,. \ea} (note that $\eta_1 = \xi_1$) and observe that, in
the new variables, the second equation of (\ref{control}) becomes
 \beeq{\label{secondeq}
 \ba{rcl} \dot \eta &=& A(\eta + M(X)\tilde
 \theta)  + \phi(\eta_1) + \Omega(\eta_1)( \theta(\varrho) + \tilde\theta) +
 H(X,\eta_1)v - M(\dot  X)
 \tilde \theta - M(X)\beta(X,\eta_1) v\\[2mm]
 &=& A\eta + [AM(X) + \Omega(\eta_1) - M(\dot  X)]\tilde \theta +
 [H(X,\eta_1)-M(X)\beta(X,\eta_1)]v + \phi(\eta_1) + \Omega(\eta_1) \theta(\varrho)\,.\ea
 }
 The third equation, instead, becomes
 \[
 \dot {\tilde \theta} = \beta(X,\eta_1) v-\mbox{dzv}_\ell\,(\tilde \theta +
 \theta(\varrho))\,.
 \]
The choices of $H(X,\eta_1), \beta(X,\eta_1)$ and of $F, G$ are
meant to simplify the terms
\[
[AM(X) + \Omega(\eta_1) - M(\dot  X)]\tilde \theta +
[H(X,\eta_1)-M(X)\beta(X,\eta_1)]v
\] in the expression (\ref{secondeq}).
First of all, note that choosing
\[
H(X,\eta_1) = M(X)\beta(X,\eta_1) + K
\]
with $K$ a constant vector, to be fixed later, the second term
becomes equal to $Kv$. The first term, instead, can be made equal
to
\[
[AM(X) + \Omega(\eta_1) - M(\dot X)]\tilde \theta= b \beta\tr
(X,\eta_1)\tilde \theta
\]
in which $b$ is a $d\times 1$ fixed vector, if $M(X)$ satisfies
\[
M(\dot X) = (A-bCA)M(X) + (I-bC)\Omega(\eta_1)
\]
and $\beta(X,\eta_1)$ satisfies \[ \beta\tr (X,\eta_1)= CAM(X) +
C\Omega(\eta_1)\,.
\]

In this way, the second equation of (\ref{control}) reduces to
\beeq{\label{MarTom} \dot \eta = A\eta + b\beta\tr
(X,\eta_1)\tilde \theta + Kv+\phi(\eta_1) + \Omega(\eta_1)
\theta(\varrho)\,. } To show that the required differential
equation for $X$ can be enforced, pick a column vector $b = {\rm
col}(1,b_2, \ldots, b_d)$. Then, bearing in mind the definition of
$M(X)$, it is easily realized that the required differential
equation holds if the matrices $F$ and $G$ have the form (see
\cite{MaTo92}) \beeq{\label{effegi} F = \qmx{-b_2 & 1 & \cdots & 0
& 0\cr \cdot & \cdot & \cdots & \cdot & \cdot\cr -b_{d-1} & 0 &
\cdots & 0 & 1\cr -b_{d} & 0 & \cdots & 0 & 0\cr}, \qquad G =
\qmx{-b_2 & 1 & \cdots & 0 & 0 & 0\cr \cdot & \cdot & \cdots &
\cdot & \cdot & \cdot\cr -b_{d-1} & 0 & \cdots & 0& 1 & 0\cr
-b_{d} & 0 & \cdots & 0 & 0 & 1\cr}.}

In summary,  the quantities $H(X,\xi_1),\beta(X,\xi_1),F,G$ which
appear in the controller (\ref{control}) are determined as
follows: $F$ and $G$ are the matrices in (\ref{effegi}), while
$\beta(X,\xi_1)$ and $H(X,\xi_1)$ are chosen as
\beeq{\label{biggamma} \beta(X,\xi_1) = [CA\qmx{0 \cr X\cr} +
C\Omega(\xi_1)]\tr } \beeq{\label{acca} H(X,\xi_1) = \qmx{0 \cr
X\cr}[CA\qmx{0 \cr X\cr} + C\Omega(\xi_1)]\tr + K\,.} The vectors
$b$  and $K$ are left undetermined, for the time being.

The controller thus defined determines a closed loop system which,
in the coordinates indicated above, can be written as
 \beeq{\label{closedloop}\ba{rcl}
 \dot \varrho &=& 0 \\[2mm]
 \dot w &=& s(\varrho,w) \\[2mm]
 \dot z &=& f_0(\varrho, w, z)+f_1(\varrho,w,z,e)e\\[2mm]
 \dot e &=& q(\varrho,w,z,e)+\eta_1+v\\[2mm]
 \dot \eta &=& A\eta + b\beta\tr (X,\eta_1)\tilde \theta + K v + \phi(\eta_1) + \Omega(\eta_1)\theta(\varrho)
 \\[2mm] \dot {\tilde \theta} &=& \beta(X,\eta_1) v -\mbox{dzv}_\ell\,(\tilde \theta + \theta(\varrho))\\[2mm]
\dot X &=& FX + G\Omega(\eta_1)\,.\ea} This system, viewed as a
system with input
 $v$ and output $e$, has relative degree 1 and a zero dynamics
 characterized by the equations
 \beeq{\label{zerodyn}\ba{rcl}
 \dot \varrho &=& 0 \\[2mm]
 \dot w &=& s(\varrho,w) \\[2mm]
 \dot z &=& f_0(\varrho, w,z) \\[2mm]
 \dot \eta &=& A\eta - K[q(\varrho,w,z,0)+\eta_1]+
 b\beta\tr (X,\eta_1)\tilde \theta + \phi(\eta_1) + \Omega(\eta_1)\theta(\varrho)
 \\[2mm] \dot {\tilde \theta} &=& -\beta(X,\eta_1) [q(\varrho,w,z,0)+\eta_1]
 -\mbox{dzv}_\ell\,(\tilde \theta + \theta(\varrho))
 \\[2mm]
 \dot X &=& FX + G\Omega(\eta_1)\,.
 \ea
 }
 Standard arguments suggest that  if the latter have convenient asymptotic properties,
 in particular possess
 a locally exponentially stable compact attractor, an additional control of the form $v=-ke$,
 (with large $k>0$) should be able to keep trajectories bounded and steer $e(t)$ to
 zero. Thus, we proceed to analyze the properties of
 (\ref{zerodyn}). The intuition that the expected asymptotic
 properties hold is based on the observation that, in view of
Assumption (iv), the fourth equation of (\ref{zerodyn}) can be
regarded as the equation of an adaptive observer, constructed
according to the methods of \cite{BaGe} and \cite{MaTo92}, for the
variable $\tau({\bf z})$, with the fifth equation providing of the
appropriate adaptation law.

\subsection{Trajectories of (\ref{zerodyn}) are bounded}

Let the initial conditions for $\varrho,w,z$ be taken in the set
${\bf Z} $, a subset of a set ${\mathcal Z}$ which by  hypothesis
is positively invariant for the (autonomous) subsystem consisting
the top three equations of (\ref{zerodyn}) and consider the change
of variables
\[
\chi = \eta - \tau(\varrho,w,z)\,.
\]
Using  (\ref{adaptimmers1}) and (\ref{adaptimmers0}) (which, it
should be borne in mind, hold only on ${\mathcal A}_0$ ) system
(\ref{zerodyn}) becomes

 \beeq{\label{zerodyn2}
 \ba{rcl}
 \dot \varrho &=& 0 \\[2mm]
 \dot w &=& s(\varrho,w) \\[2mm]
 \dot z &=&  f_0(\varrho,w,z) \\[2mm]
 \dot \chi &=& (A-KC)\chi+ b\beta\tr (X,\chi_1+\tau_1)\tilde \theta +
 \Delta(\chi_1,\tau_1,\theta) + \varphi(\varrho, w,z)\\[2mm]
 \dot {\tilde \theta} &=& -\beta(X,\chi_1+\tau_1)\chi_1 -\mbox{dzv}_\ell\,(\tilde \theta +
 \theta(\varrho))\\[2mm]
 \dot X &=& FX + G\Omega(\chi_1+\tau_1)\,,\ea} in which
 \[
 \Delta(\chi_1,\tau_1,\theta) = \phi(\chi_1+\tau_1) - \phi(\chi_1)
 + [\Omega(\chi_1+\tau_1) - \Omega(\chi_1)]\theta(\varrho)
 \]
 is a term which vanishes at $\chi_1=0$ and
 \beeq{\label{defe}\ba{rcl}
 \varphi(\varrho, w,z) &=& K(c(\varrho, w,z)-\tau_1(\varrho, w,z)) + A \tau(\varrho, w,z)) + \phi(\tau_1(\varrho, w,z))
 + \Omega(\tau_1(\varrho, w,z)) \theta(\varrho) \\[2mm] &&\dst \qquad-
 \,{\partial \tau \over \partial {z}} {
 f}_0(\varrho, w, z)
 - {\partial \tau \over \partial {w}} {
 s}(\varrho, w)\ea
 }
 is a term vanishing on ${\mathcal A}_0$.
 In particular note that, since $\phi(\cdot)$ and $\Omega(\cdot)$ can be taken to be
globally Lipschitz and $\theta$ ranges over a compact set, there
exists a number $L$ such that
\[
|\Delta(\chi_1,\tau_1,\theta)| \le
L_\phi\,|\chi_1|+L_\Omega|\chi_1|\,|\theta|\le L|\chi_1|\] for all
$\chi_1,\tau_1,\theta$.

Mimicking a construction used in \cite{MaTo92} for the design of
adaptive observers, the vectors $b_i$'s and $K$ are chosen as
follows. The $b_i$'s be such that the polynomial
 \beeq{ \label{polync}
 p(\lambda) = \lambda^{d-1} + b_2\lambda^{d-2} + \cdots +
 b_{d-1}\lambda + b_d
 }
 has $d-1$ distinct roots with negative real part, in which case
 the matrix $F$ in the bottom equation of
(\ref{zerodyn2}) is Hurwitz (and has distinct eigenvalues). On the
other hand, $K$ is chosen as \beeq{\label{kappa} K = Ab + \lambda
b } in which $\lambda
>0$. These choices are such that the following result holds.

\begin{lemma}\label{LM1} Suppose assumptions (i), (ii), (iv) hold. There is a number $\lambda^\ast$ such
that, if $\lambda \ge \lambda^\ast$, all trajectories of
(\ref{zerodyn2}) are bounded.
\end{lemma}

\begin{proof}
Because of assumption (ii), $(\varrho, w(t), z(t)) \in {\mathcal
Z}$ for all $t \geq 0$, where ${\mathcal Z}$ is a compact set.
Form this, standard arguments can be used to show  that
trajectories of (\ref{zerodyn2}) are defined for all $t\ge 0$ and,
consequently, also $X(t)$ is bounded (recall that
$|\Omega(\cdot)|$ has compact support). To prove the Lemma it
remains to show that also $\chi(t)$ and $\tilde \theta(t)$ are
bounded. To this end, let $\chi$ be partitioned as $\chi = {\rm
col}(\chi_1,\chi_2)$, in which $\chi_2$ is a $(d-1)\times 1$
vector and replace $\chi_2$ by
 \[
 \zeta = \hat b \chi_1 + \chi_2\,,
 \]
 where $\hat b = -{\rm col}(b_2, b_3, \ldots, b_d)\,.$
 In this way, the fourth and fifth equations of
 (\ref{zerodyn2}) are changed to
 \beeq{\label{MarTom2}\ba{rcl}
 \dot \chi_1 &=& -\lambda \chi_1 + \hat c\zeta +
 \beta\tr (X,\chi_1+\tau_1)\tilde \theta + C\Delta(\chi_1,\tau_1,\theta) + C \varphi(\varrho, w, z)\\[2mm]
 \dot \zeta &=& F\zeta + \qmx{\hat b & I \cr}\Delta(\chi_1,\tau_1,\theta)
 +\qmx{\hat b & I \cr} C \varphi(\varrho, w, z)\\[2mm]
 \dot {\tilde \theta} &=& -\beta (X,\chi_1+\tau_1)\chi_1 - \mbox{dzv}_\ell(\tilde \theta + \theta(\varrho))\,.
 \ea}

 Choose now for (\ref{MarTom2}) the Lyapunov function
 \beeq{\label{Vdef}
 V(\chi_1,\zeta,\tilde \theta) = \chi_1^2 + \zeta\tr P \zeta +
 \tilde \theta\tr \tilde \theta\,,
 }
 in which $P$ is the positive definite solution of $PF+F\tr P = -
 I$ and obtain, after some simple algebra,
 \beeq{\label{dotv.1}\ba{rcl} \dot V
 &\le& -2\lambda \chi_1^2 - |\zeta|^2 - 2 \tilde \theta\tr
 \mbox{dzv}_\ell(\tilde \theta + \theta(\varrho)) + L_1|\chi_1|^2 +
 L_2|\chi_1|\,|\zeta|\\[2mm]
 && +L_3|\chi_1||\varphi(\varrho, w, z)| + L_4 |\zeta||\varphi(\varrho, w,
 z)|
 \ea}
 for some $L_i>0$, $i=1\ldots,4$.
 Bearing in mind the fact that $\varphi(\varrho, w, z)$ is
 bounded by some fixed number $\bar \varphi>0$ and
completing the squares, one obtains
 \beeq{\label{Vdot1}
 \dot V \le -(2\lambda -L_1 + {\dst 1 \over \dst 2}L_2^2) \chi_1^2  - {\dst 1 \over 2}|\zeta|^2
 - 2 \tilde \theta\tr \mbox{dzv}_\ell(\tilde \theta + \theta(\varrho))
 +L_3\bar |\chi_1|\varphi + L_4\bar  |\zeta| \varphi\,.
 }
A property of the function (\ref{deftas}), in view of the choice
of $\ell$, is that \beeq{\label{tas2}
  \tilde \theta\tr \mbox{dzv}_\ell(\tilde \theta +
 \theta(\varrho)) \geq 0 \qquad \mbox{for all } \tilde \theta \in
 \Real^q\quad \mbox{and } \quad \varrho \in P\,.
 }
 Moreover, for any $\delta >  \sqrt{q}(2 \ell+1)$ there is a
 number $c_1>0$ such that
 \beeq{\label{tas1}
 |\tilde \theta| \geq \delta \qquad \Rightarrow \qquad
 2 \tilde \theta\tr \mbox{dzv}_\ell(\tilde \theta + \theta(\varrho)) \geq
 c_1 |\tilde \theta|^2 \qquad \mbox{for all } \tilde \theta \in
 \Real^q\quad \mbox{and } \quad \varrho \in P\,.
 }
If $\lambda$ is large enough, inequality (\ref{Vdot1}), in view of
property (\ref{tas1}), yields
 \[
 |\tilde \theta| \geq \delta \qquad \Rightarrow \qquad \dot V \leq
 -c_2  |(\chi_1, \zeta, \tilde \theta)|^2  + c_3 |(\chi_1, \zeta, \tilde
 \theta)|\bar \varphi
 \]
 for suitable $c_2 >0, c_3>0$. This, in turn, yields
 \beeq{\label{Re1}
  |\tilde \theta| \geq \delta \quad \mbox{and} \quad
  |(\chi_1, \zeta, \tilde \theta)| > {c_3 \over c_2} \, \bar \varphi
 \qquad \Rightarrow \qquad \dot V <0\,.
 }
Property (\ref{tas2}), on the other hand, yields
 \[
 \dot V \leq - c_2 |(\chi_1, \zeta)|^2 + c_3 |(\chi_1, \zeta)| \,
 \bar \varphi
 \]
which, in turn, yields
 \beeq{\label{Re2}
 |(\chi_1, \zeta)| > {c_3 \over c_2} \, \bar \varphi
 \qquad \Rightarrow \qquad \dot V <0\,.
 }
At this point, it is easy to conclude that (\ref{Re1}) and
(\ref{Re2})
  imply
 \[
 |(\chi_1, \zeta, \tilde \theta)| > \sqrt{\delta^2 + \Bigl({c_3 \over c_2} \,\bar \varphi\Bigr)^2}
  \qquad \Rightarrow \qquad
 \dot V <0\,.
 \]
 Bearing in mind the fact that $V(\chi_1,\zeta,\tilde \theta)$ is a quadratic form,
  the claim follows by standard arguments.
 $\triangleleft$
\end{proof}

\subsection{The limit set of (\ref{zerodyn}) and its properties}\label{sec3}

Let the  initial conditions $\eta(0),\tilde \theta (0), X(0)$ of
(\ref{zerodyn}) be taken in fixed compact sets ${\bf H},\,
\Theta,\, {\bf X}$. From Lemma \ref{LM1} we can claim that, if
$\lambda$ is large enough, the positive orbit of the set
\[
{\bf B} = P\times W\times Z \times {\bf H} \times \Theta \times
{\bf X}
\]
under the flow of (\ref{zerodyn}) is bounded. As a consequence
$\omega({\bf B})$, the $\omega$-limit set of ${\bf B}$ under the
flow of (\ref{zerodyn}), is a non-empty, compact and invariant
set, which uniformly attracts  all trajectories of (\ref{zerodyn})
with initial conditions in ${\bf B}$. In what follows we render
the structure of $\omega({\bf B})$ explicit. To this end, we need
an extra hypothesis, which -- as any well-educated reader in
adaptive control is expecting -- plays in the present context a
role completely similar to the role of the classical hypothesis of
{\em persistence of excitation}.

Consider the equivalent system (\ref{zerodyn2}) and rewrite the
three top equations as in (\ref{zerobold}) (consistently rewrite
$\varphi(\varrho, w,z)$ as $\varphi({\bf z})$ and
$\tau(z,w,\varrho)$ as $\tau({\bf z})$). Because of the special
triangular structure of (\ref{zerodyn2}), it can be observed that
if $({\bf z},\chi,\tilde \theta, X)$ is a point of $\omega({\bf
B})$, necessarily ${\bf z}$ is a point in the $\omega$-limit set
of ${\bf Z}$ under the flow of (\ref{zerobold}), that is, ${\bf
z}$ is a point of ${\mathcal A}_0$. This implies that on
$\omega({\bf B})$ we have $\varphi({\bf z})=0$ and thus system
(\ref{zerodyn2}) simplifies as
 \beeq{\label{zerodyn.2}
 \ba{rcl}\dot
 {\bf z} &=& {\bf f}_0({\bf z})\\[2mm]\dot \chi &=& (A-KC)\chi + b\beta\tr (X,\chi_1+\tau_1)\tilde
 \theta +
 \Delta(\chi_1,\tau_1,\theta) \\[2mm] \dot {\tilde \theta} &=& -\beta(X,\chi_1+\tau_1)
 \chi_1 - \mbox{dzv}_\ell(\tilde \theta + \theta(\varrho))\\[2mm]
 \dot X &=& FX + G\Omega(\chi_1+\tau_1)\,.\ea
 }

 It will be shown now that, if the announced additional
 hypothesis of persistency of excitation holds,
 points of $\omega({\bf B})$ have necessarily $\chi=0$, $\tilde
 \theta=0$. To introduce the hypothesis in question
 we observe first of all the following interesting feature.

\begin{lemma}\label{ssX}
The graph of the map
\[
\ba{rcccl}
 \sigma &:& {\mathcal  A}_0 &\to& \Real^{(d-1)\times q}\\
 && {\bf z} & \mapsto & \dst X=\int_{-\infty}^0
 e^{-Fs}G\Omega(\tau_1({\bf z}(s,{\bf z}))) ds
 \ea
\]
 is invariant for \beeq{\label{zerodyn3} \ba{rcl} \dot
 {\bf z} &=& {\bf f}_0({\bf z})\\[2mm]
 \dot X &=& FX + G\Omega(\tau_1({\bf z}))\,.\ea
 }
\end{lemma}

\begin{proof}
 Let ${\bf z}(t,{\bf z}_0)$ denote the solution of (\ref{zerobold}) passing
through ${\bf z}_0$ at time $t=0$ and note that, if ${\bf z}_0 \in
{\mathcal A}_0$, then ${\bf z}(t,{\bf z}_0) \in {\mathcal A}_0$
for all $t$. Since $F$ is a Hurwitz matrix, the map
$\sigma(\cdot)$ is well defined. A simple calculation shows that
\[
 \sigma({\bf z}(t,{\bf z}_0)) =e^{F t} \sigma({\bf z}_0) +
 \int_0^t e^{F(t- s)} G \Omega(\tau_1({\bf z}(s,{\bf z}_0)))
 ds\,,
\]
from which the result follows. $\triangleleft$ \end{proof}

The new assumption is the following one.

\bigskip\noindent
{\em Assumption (v)}\,: Consider the map $\gamma: {\mathcal A}_0
\to \Real^{q\times 1}$ defined as \[ \gamma:  {\bf z}
\;\;\;\mapsto\;\;\; \beta(\sigma({\bf z}),\tau_1({\bf z}))\] It is
assumed that for any initial condition ${\bf z}_0 \in {\mathcal
A}_0$
 the identity
\[
c\tr \gamma({\bf z}(t,{\bf z}_0)) = 0, \qquad \mbox{for all $t\in
\Real$}\] implies $c=0$. $\triangleleft$

\bigskip
Under this hypothesis, the set $\omega({\bf B})$ assumes a very
simple structure. As a matter of fact, the following result holds.

\begin{proposition}\label{PR1} Under the assumptions (i),(ii),(iv) and (v)
the set $\omega({\bf B})$ is the graph of a continuous map defined
on ${\mathcal A}_0$. Any point of $\omega({\bf B})$ is a point
$({\bf z},\eta,\tilde \theta,X)$ in which ${\bf z}\in {\mathcal
A}_0$ and
\[\eta=\tau({\bf z}), \qquad
\tilde \theta = 0, \qquad X = \sigma({\bf z})\,.\] If also
assumption (iii) holds, then  $\omega({\bf B})$ is locally
exponentially attractive for (\ref{zerodyn}).\end{proposition}

\begin{proof}  By contradiction, suppose
a point ${\bf p}=({\bf z},\chi_0,\tilde \theta_0, X)$ with either
$\chi_0\ne 0$ or $\tilde \theta_0\ne 0$ is in $\omega({\bf B})$.
Since $\omega({\bf B})$ is compact and invariant, the backward
trajectory of (\ref{zerodyn.2}) starting at this point is bounded.
Along this trajectory, the Lyapunov function
  (\ref{Vdef}) satisfies $V \le C$ for all $t\le 0$, for some $C>0$.
 Moreover, since $\varphi({\bf z})=0$ on $\omega({\bf B})$, the same
 computations indicated in the proof of Lemma \ref{LM1} show that

 \[
 \dot V(t) \leq -(2 \lambda - L_1 - {1 \over 2} L_2^2)|\chi_1(t)|^2
 -{1 \over 2} |\zeta(t)|^2 - 2\tilde \theta(t) \mbox{dzv}_{\ell}(\tilde \theta(t) +
 \theta(\varrho))\,.
 \]
  Using property (\ref{tas2})
  it is seen that,  if $\lambda \geq \lambda^\star$, $V(t)$ is
 non-increasing along trajectories.  As consequence, since $V(t)$ is bounded, there
 exists a finite number $V_\alpha$ such that $\lim_{t\to -\infty}V(t)=V_\alpha\,.$ The trajectory in question
is attracted, in backward time, by its own $\alpha$-limit set
$\alpha({\bf p})$, which, as it is well known, is nonempty,
compact and invariant. Moreover, by definition, the function
$V(\chi_1,\zeta,\tilde \theta)$ has the same value $V_\alpha$ at
any point of $\alpha({\bf p})$. Proceeding as in the classical
proof of LaSalle's invariance principle, pick an initial condition
$\hat {\bf p}$ in the set $\alpha( {\bf p})$ and consider the
corresponding trajectory of (\ref{zerodyn.2}). Along such
trajectory, $V(t)$ is constantly equal to $V_\alpha$ and hence
\[
\chi_1(t) = 0, \qquad \zeta(t)=0, \qquad
 \mbox{dzv}_\ell(\tilde \theta(t) + \theta(\varrho))=0\qquad
\mbox{for all $t\in \Real$.}\] Entering these constraints in
(\ref{zerodyn.2}), and observing that the vector $b$ is nonzero,
it is seen that necessarily
\[\ba{l}
\tilde \theta \tr \beta  =0\\[2mm]
\dot{\tilde \theta} = 0\\[2mm]
\dot X = FX + G\Omega(\tau_1({\bf z}))\,.\ea\] It is seen from
this that $\tilde \theta(t)$ is a constant, say $\tilde
\theta^\ast$, along such trajectory, while $X(t)$ is a solution of
(\ref{zerodyn3}). Since $F$ is Hurwitz and has distinct
eigenvalues, $X(t)$ can be bounded for $t\le 0$ only if $X(0) =
\sigma({\bf z}(0))$, where $\sigma(\cdot)$ is the map introduced
in Lemma \ref{ssX}, in which case $X(t) = \sigma({\bf z}(t))$.
Since $X(t)$ has to be bounded because $\alpha({\bf p})$ is
compact, it follows that $X(t)$ is necessarily equal to
$\sigma({\bf z}(t))$. This being the case,
 bearing in mind  the expression of $\beta(\cdot,\cdot)$ and the definition of the
map $\gamma(\cdot)$, the first condition shows that necessarily
\[
 ( \tilde\theta^\ast)\tr \gamma({\bf z}(t))= 0, \qquad \mbox{for all $t\in
\Real$}.\] Thus, from Assumption (v), it is concluded that $\tilde
\theta^\ast =0$. As a consequence $(\chi_1,\zeta,\tilde
\theta)=(0,0,0)$ at any point of $\alpha( {\bf p})$, and $V_\alpha
=0$. But this is a contradiction, because $V(t)$ is non-increasing
along trajectories and $V(0)$ is strictly positive, if either
$\chi_0 \ne 0$ or $\tilde \theta_0\ne 0$.

The proof that $\omega({\bf B})$ is locally exponentially stable
is a consequence of Assumption (iii). It requires appropriate
modifications of arguments used in similar instances in the proof
of convergence of parameter estimates in various adaptive control
schemes, such as those presented in various Chapters of
\cite{Khalil}, and is not included here for reasons of space.
 Details can be found in \cite{DMISIAM}.
$\triangleleft$
\end{proof}

\section{Adaptive output regulation}

We return now to the closed loop system obtained from the
interconnection of (\ref{pla1}), (\ref{exo1}) and (\ref{control}).
 As observed, this system, viewed as a system with input
$v$ and output $e$ has relative degree 1. To put it in ``normal
form" we use, instead of (\ref{newvar}), the change of variables
\beeq{\label{newvar2}\ba{rcl}\tilde \theta &=& \hat \theta -
\theta(\rho) - \dst\int_0^e\beta(X,\xi_1-CKe + CKs)ds
\\[2mm]
\eta &=& \xi - M[\hat \theta - \theta(\rho)] - Ke\,.\ea}

This, after some simple algebra and some obvious rearrangement of
terms, yields a system of the form
\beeq{\label{closedloop2}\ba{rcl}
 \dot \varrho &=& 0 \\[2mm]
 \dot w &=& s(\varrho,w) \\[2mm]
 \dot z &=&  f_0(\varrho,w,z)+f_1(\varrho, w,z,e)e \\[2mm]
 \dot \eta &=& A\eta +
 b\beta\tr \tilde \theta(\varrho) - K[q(\varrho,w,z,0)+\eta_1] +
 \phi(\eta_1) + \Omega(\eta_1)\theta +
 \delta_1(\varrho,w,z,e,X,\eta_1)\,e\\[2mm]
\dot {\tilde \theta} &=& -\beta [q(\varrho,w,z,0)+\eta_1] +
\delta_2(\varrho,w,z,e,X,\eta_1)\,e - \mbox{dzv}_\ell(\tilde
\theta + \theta(\varrho))\\[2mm]
\dot X &=& FX + G\Omega(\eta_1)+\delta_3(\eta_1,e)\,e\\
[2mm]\dot e &=&
-[q(\varrho,w,z,0)+\eta_1]+\vartheta(\varrho,w,z,e)e+v\,.\ea} in
which $\delta_1(\cdot)$, $\delta_2(\cdot)$, $\delta_3(\cdot)$ and
$\vartheta(\cdot)$ are continuously differentiable functions of
their arguments. The notation $\beta$ stands for
$\beta(X,\eta_1+CKe)$.

This system can be put in a  much more compact form by setting
 \[\ba{rcl} {\bf x} &=& {\rm col}(\varrho,w,z,\eta,\tilde \theta, X_1,\ldots X_q)\,, \ea\]
(where $X_i$ denotes the $i$-th column of $X$) which yields a
system of the form
 \beeq{\label{closbold}
 \ba{rcl}
 \dot {\bf x} &=& {\bf F}_0({\bf x})+ {\bf F}_1({\bf
 x},e)e
 \\[2mm]\dot e &=& {\bf h}({\bf x})+ {\bf k}({\bf x}, e)e+v\,.
 \ea}
 In this notation, the ``subsystem" $\dot {\bf x} = {\bf F}_0({\bf
 x})$
 is precisely system (\ref{zerodyn}),
while the function ${\bf h}({\bf x})$, which is precisely the
quantity $-[q(\varrho,w,z,0 )+\eta_1]$ in (\ref{closedloop2}),
vanishes on the set $\omega({\bf B})$. Having realized this, it is
not difficult to claim that the controller (\ref{control})
completed with
 \beeq{\label{vgx}
 v = - k e\,,
 }
 if $k$ is sufficiently large, keeps trajectories bounded and
 steers $e(t)$ to 0.

\begin{proposition} \label{Main} Consider system (\ref{pla1}) with exosystem (\ref{exo1}).
Let  ${\bf Z}, E$ be fixed compact sets of initial conditions, for
which the assumptions (i)-(iv) indicated in section \ref{SecPA}
are supposed to hold. Suppose, in addition, that assumption (v)
introduced in section \ref{sec3} holds. Consider the controller
(\ref{control}) completed with (\ref{vgx}) and initial conditions
in a fixed compact set ${\bf K}$. Then, there exists a number
$k^\ast>0$ such that if $k\ge k^\ast$ the positive orbit of ${\bf
Z} \times  E \times {\bf K}$ in the closed loop system is bounded
and $e(t)\to 0$ as $t\to \infty$.
\end{proposition}

\begin{proof}
The closed-loop system \beeq{\label{closbold2}
 \ba{rcl}
 \dot {\bf x} &=& {\bf F}_0({\bf x})+ {\bf F}_1({\bf
 x},e)e
 \\[2mm]\dot e &=& {\bf h}({\bf x})+ {\bf k}({\bf x}, e)e-ke\,.
 \ea}
can be viewed as interconnection of two subsystems, one with state
${\bf x}$ and input $e$, the other with state $e$ and input ${\bf
x}$. As shown in \cite{BIP}, the upper subsystem is input-to-state
stable (relative to the set $\omega(B)$), with a gain function
which is linear if the assumptions (i) through (v) hold. Then, an
easy extension (to the case of systems which are input-to-state
stable relative to compact attractors) of the small-gain theorem
of \cite{JTP} can be invoked to show that, if $k$ is large enough,
the results of the Proposition hold. Further details on the proof
of this ``high-gain" stabilizability property can be found in
\cite{BIP}.
\end{proof}$\triangleleft$

Before concluding the paper, we wish to point out the fact that
the theory developed so far lends itself to deal with more general
situations in which the equations of the (augmented) {\em
exosystem} are {\em affected by regulated variable} $e$.
Specifically, suppose that $\varrho$ and $w$, instead of
(\ref{exo1}), are generated by a system of the more general form
\beeq{\label{exo2}\ba{rcl}
 \dot \varrho&=& s_\varrho(\varrho,w,z,e)e\\
 \dot w &=& s(\varrho,w)+s_w(\varrho,w,z,e)e \,.
 \ea}

In this case, in fact, the closed loop system obtained from the
interconnection of (\ref{pla1}), (\ref{exo2}) and (\ref{control})
can still be put in the form (\ref{closbold}), in which ${\bf
f}_0({\bf z})$ has exactly the form (\ref{zerobold}), while
\[
{\bf f}_1({\bf z},e) = \qmx{s_\varrho(\varrho,w,z,e)\cr
s_w(\varrho,w,z,e)\cr f_1(\varrho,w,z,e)\cr}\,.
\]
Thus, the result expressed by Proposition \ref{Main} continues to
hold. Models of this kind occur, for instance, in the problem of
self-compensation of mechanical and/or electrical asymmetries in a
rotating electrical drives (see \cite{BIMP}).

\section{Conclusions}

This paper has discussed the design of nonlinear internal models
in the problem of adaptive output regulation. It has been shown
how the theory of adaptive observers can be successfully used to
deal with complex situations, not covered by existing results, in
which the desired steady state control input is generated by a
possibly nonlinear system and depends on constant uncertain
parameters. The result is framed in the general
``non-equilibrium'' theory proposed in \cite{BI03}, thus allowing
for controlled plant with not necessarily stable zero dynamics.

\newpage
\begin{center} {\bf Acknowledgements}
\end{center}
\medskip
The authors wish to thank Christopher I. Byrnes and Laurent Praly
for fruitful discussion and helpful suggestions during the
preparation of the paper.

\end{document}